\documentclass{amsart}
\usepackage{graphicx, amssymb}
\usepackage[all,xdvi]{xy}
\newtheorem*{theorem*}{Theorem}
\newtheorem{theorem}{Theorem}[section]
\newtheorem{lemma}[theorem]{Lemma}

\newtheorem{corollary}[theorem]{Corollary}

\theoremstyle{definition}
\newtheorem{definition}[theorem]{Definition}

\newtheorem{example}[theorem]{Example}

\theoremstyle{remark}
\newtheorem{remark}[theorem]{Remark}

\numberwithin{equation}{section}
\newcommand{\firef}[1]{Figure~{\rm\ref{#1}}}
\newcommand{\thref}[1]{Theorem~{\rm\ref{#1}}}

\newcommand{\leref}[1]{Lemma~{\rm\ref{#1}}}
\newcommand{\coref}[1]{Corollary~{\rm\ref{#1}}}

\newcommand{\deref}[1]{Definition~{\rm\ref{#1}}}
\newcommand{\exref}[1]{Example~{\rm\ref{#1}}}

\newcommand{\seref}[1]{Section~{\rm\ref{#1}}}

\newcommand{\fig}[1]
{\raisebox{-0.5\height}%
{\includegraphics{#1}}}
\newcommand{\st}{\; | \;}                               
\newcommand{\ttt}{\otimes}                              

\newcommand{\<}{\langle}
\renewcommand{\>}{\rangle}
\newcommand{\isoto}{\xrightarrow{\sim}}       
\newcommand{\xxto}{\xrightarrow}              

\renewcommand{\i}{{\mathrm{i}}}   
\newcommand{\CC}{\mathbb{C}}       
\newcommand{\Z}{\mathbb{Z}}       
\newcommand{\PP}{\mathbb{P}^1}    
\newcommand{\OO}{\mathcal{O}}     
\newcommand{\F}{\mathcal{F}}      
\newcommand{\G}{\mathcal{G}}      
\newcommand{\TT}{\mathcal{T}}      
\newcommand{\D}{\mathcal{D}}      
\newcommand{\C}{\mathcal{C}}      

\newcommand{\al}{\alpha}

\newcommand{\Ga}{\Gamma}
\newcommand{\de}{\delta}
\newcommand{\De}{\Delta}

\newcommand{\la}{\lambda}

\newcommand{\ph}{\varphi}
\newcommand{\Ph}{\Phi} 
\newcommand{\si}{\sigma}

\newcommand{\eps}{\varepsilon}

\newcommand{\g}{\mathfrak{g}}                  
\newcommand{\n}{\mathfrak{n}} 
 
\newcommand{\slt}{\mathfrak{sl}_2}             

\newcommand{\Qhat}{{\widehat{Q}}}
\newcommand{\Gcoh}{{Coh_G(\PP)}}
\newcommand{\Gbar}{{\bar G}}                  
\newcommand{\Gbarcoh}{{Coh_{\Gbar}(\PP)}}
\newcommand{\Gqcoh}{{Qcoh_G(\PP)}}
\newcommand{\Coh}{Coh(\PP)}
\newcommand{\Qcoh}{Qcoh(\PP)}
\newcommand{\SU}{\mathrm{SU(2)}}

\DeclareMathOperator{\Rep}{Rep}
\DeclareMathOperator{\Ind}{Ind}
\DeclareMathOperator{\Irr}{Irr}
\DeclareMathOperator{\id}{id}

\DeclareMathOperator{\Hom}{Hom}
\DeclareMathOperator{\RHom}{RHom}
\newcommand{\GHom}{\Hom_G}
\newcommand{\GExt}{\Ext_G}
\DeclareMathOperator{\Ext}{Ext}

\DeclareMathOperator{\rk}{rk}
\DeclareMathOperator{\Stab}{Stab}
\DeclareMathOperator{\Path}{Path_{\Qhat}}

\begin{document}

\title{McKay correspondence and equivariant sheaves on $\PP$}

\author{Alexander Kirillov, Jr.}
   \address{Department of Mathematics, SUNY at Stony Brook, 
            Stony Brook, NY 11794, USA}
    \email{kirillov@math.sunysb.edu}
    \urladdr{http://www.math.sunysb.edu/\textasciitilde kirillov/}
\dedicatory{Dedicated to my father on the occasion of his 70th birthday}

\begin{abstract}
Let $G$ be a finite subgroup in $\SU$, and $Q$ the corresponding
affine Dynkin diagram. In this paper, we review the relation between
the categories  of $G$-equivariant sheaves on $\PP$ and $\Rep
Q_h$, where $h$ is an orientation of $Q$, constructing an explicit
equivalence  of corresponding derived categories.  
\end{abstract}
\maketitle
\section*{Introduction}
Let $G\subset \SU$ be a finite subgroup. According to McKay correspondence,
such a subgroup gives rise to a graph $Q$ which turns out to be an affine
Dynkin diagram of ADE type. Let $\g$ be the corresponding affine Lie algebra.

There are several approaches allowing one to  construct $\g$ from $G$: 
\begin{enumerate}
\item There is a geometric construction, due to Nakajima and Lusztig,  of $\g$
  and its representations in  terms of the quiver varieties associated to graph
  $Q$. These varieties are closely related to the moduli spaces of instantons
  on the resolution of singularities $\widetilde{\CC^2/G}$. 

\item There is an algebraic construction, due to Ringel \cite{ringel2}, which
  allows one to  get the universal enveloping algebra $U\n_+$ of the positive
  part of $\g$ as the Hall algebra of the category $\Rep Q_h$, where $h$ is an
  orientation of $Q$ and $Q_h$ is the corresponding quiver. It was shown by 
  Peng and Xiao \cite{peng-xiao} that replacing $\Rep Q_h$ by the
  quotient  $R=D^b(\Rep Q_h)/T^2$ of the corresponding derived category, we
  can get a description of all of $\g$.  Different choices of orientation
  give rise to equivalent derived categories  $R$, with equivalences given
  by  Bernstein-Gelfand-Ponomarev reflection functors; in
  terms of $\g$, these functors correspond to the action of the braid group
  of the corresponding Weyl group.
  
  \end{enumerate}
Recently, a third approach was suggested by Ocneanu \cite{ocneanu}
(unpublished), in a closely related setting of a subgroup in quantum
$\SU$, for  $q$ being a root of unity. His approach is based on
studying essential paths in $\Qhat=Q\times \Z/h\Z$, where $h$ is the
order of the root of unity. This approach is purely combinatorial: all
constructions are done using this finite graph and vector spaces of
essential paths between points in this graph, without involving any
categories at all. 

This paper grew out of the author's attempt to understand Ocneanu's
construction and in particular, find the appropriate
categorical interpretation of  his combinatorial constructions;
however, for simplicity we  do it for subgroups in classical $\SU$,
leaving the analysis of the subgroups in quantum $\SU$ for future
papers. We show that Ocneanu's essential paths in $\Qhat$ have a
natural interpretation in terms of the category  $\Gcoh$ of
$G$--equivariant $\OO$-modules on $\PP$ (or, rather,
its ``even part'' $\C=\Gcoh_0$, for certain natural $\Z_2$ grading on
$\Gcoh$). This also provides a relation with Ringel--Lusztig
construction: for (almost) any choice of
orientation $h$ of $Q$, we construct an equivalence of triangulated
categories 
$$
R\Ph_h\colon D^b( \Gcoh_0)\simeq D^b(\Rep Q_h)
$$
These equivalences agree with the equivalences of $D^b(\Rep Q_h)$ for
different choices of $h$ given by reflection functors. As a corollary, we
see that the Grothendieck group $L=K(\C)$ is an affine root lattice, and
the set $\Delta$ of classes of indecomposable modules is an affine root
system. 

This construction of the affine root system via equivariant
sheaves has a number of remarkable properties,  namely:
\begin{enumerate}
  \item This does not require a choice of orientation of
    $Q$ (unlike   the category $D^b(\Rep Q_h)$, where we first choose
    an  orientation and then   prove that the resulting
    derived category is independent of orientation).
  \item The indecomposable objects in the category $\Gcoh$ can be
    explicitly described. Namely, they  are the
    sheaves $\OO(n)\ttt X_i$, where $X_i$ are irreducible
    representations of $G$ (these sheaves and their translations
    correspond to real roots of $\g$) and
    torsion sheaves, whose support is a $G$--orbit in $\PP$ (in
    particular, torsion sheaves whose support is an orbit of a
    generic  point correspond to imaginary roots of $\g$).
  \item This construction of the affine root system does not give 
    a natural  polarization into negative and positive roots.
    Instead, it gives a canonical Coxeter element in the corresponding
    affine Weyl group, which is given by $C[\F]=[\F\ttt \OO(-2)]$
    (this  also corresponds to the  Auslander--Reiten functor
    $\tau$).
 
  \item This construction gives a bijection of the vertices of the
    affine Dynkin diagram $Q$ and (some of) the $C$-orbits in
    $\Delta$ (as opposed to the quiver construction, where vertices
      of $Q$ are in bijection with the simple positive roots). 
          
\end{enumerate}

For $G=\{1\}$ (which is essentially equivalent to $G=\{\pm 1\}$),
these results were first obtained in the paper \cite{baumann-kassel},
where it was shown that the corresponding Hall algebra contains the
subalgebra isomorphic to $U\widehat{\slt}$. This  paper in turn was
inspired by an  earlier paper of Kapranov \cite{kapranov}.

It should be noted that many of the results obtained here have already
been proved in other ways. Most importantly, it had been proved by
Lenzing \cite{geigle-lenzing,lenzing} that the derived category  of
equivariant sheaves on $\PP$ is equivalent to the derived category of
representation of the corresponding quiver; this result has been used by
Schiffmann \cite{schiffmann,schiffmann2} for construction of a subalgebra
in the corresponding affine Lie algebra via Hall algebra of the category of
equivariant sheaves. However, the construction of
equivalence in \cite{lenzing} is different than the one suggested here. The
primary difference is that in Lenzing's construction, the Dynkin diagram
$Q$ is constructed as the star diagram, with lengths of branches determined
by the branching points of the cover $\PP\to X=\PP/G$, and he uses
a standard orientation of this diagram. In the construction presented
here, the diagram is $Q$ is defined in a more standard way, using the set
$I$ of irreducible representations of $G$; more importantly, we construct
not a single equivalence but a collection of equivalences, one for  each
admissible orientation. In Lenzing's construction, torsion sheaves
naturally play a major role; in our construction, we concentrate on
the study of locally free sheaves. 

Lenzing's results   apply not only to $\PP/G$ but to a much larger
class of ``non-commutative projective curves''. However, the downside of
this is that the language he uses is rather technical, making his papers
somewhat hard for non-experts. For this reason, we had chosen to  give
independent proofs of some of the results, thus saving the reader the
necessity of learning the language of non-commutative  curves. Of course,
we tried to clearly mark the results which had already been known.


{\bf Acknowledgments.} The author would like to thank P.~Etingof,
V.~Ostrik, and O.~Schiffmann for fruitful discussions. Special thanks to
A.~Ocneanu, whose talk inspired this paper. 

\section{Basic setup}\label{s:basic}
Throughout the paper, we work over the base field $\CC$ of complex
numbers. $G$ is a finite subgroup in $\SU$; for simplicity, we assume
that $\pm I\in G$ and denote  $\Gbar=G/\{\pm I\}$ (this excludes
$G=\Z_n$, $n$ odd; this case can also be included, but some of the
constructions of this paper will require minor changes). We denote by
$V$ the standard 2-dimensional space, considered as a tautological
representation of $\SU$ (and thus of $G$) and 
$$
V_k=S^k V
$$
is the $k+1$--dimensional irreducible representation of $\SU$ (and
thus a representation, not necessarily irreducible, of $G$). 

We denote by $\Rep G$ the category of finite-dimensional
representations of $G$ and by $I=\Irr G$ the set of isomorphism
classes of simple representations; for $i\in I$, we denote by $X_i$
the corresponding representation of $G$. Since $G\supset \{\pm I\}$,
the category $\Rep G$ is naturally $\Z_2$--graded:
\begin{equation}\label{e:even-odd}
\Rep G=\Rep_0 G\oplus \Rep_1 G,\quad 
\Rep_p G=\{V\in \Rep G\st (-I)|_V=(-1)^p\id \}
\end{equation}
For homogeneous object $X$, we define its ``parity'' $p(X)\in \Z_2$
by  
\begin{equation}\label{e:parity}
p(X)=p\text{ if } X\in \Rep_p G
\end{equation}
in particular, $p(V)=1$. We will also define, for $i\in I$, its
parity $p(i)=p(X_i)$. This  gives a decomposition 
\begin{equation}\label{e:decomposition_I}
I=I_0 \sqcup I_1.
\end{equation}

We define the graph $Q$ with the set of vertices $V(Q)=I$ and
for every two vertices $i,j$, the number of edges connecting them is
defined by 
$$
n(i,j)=\dim \Hom_G(X_i,X_j\ttt V)
$$

Note that decomposition \eqref{e:decomposition_I} shows that this
graph is bipartite: $V(Q)=V_0(Q)\sqcup V_1(Q)$, and all edges connect
vertices of different parities.

It is well-known that  one can construct an isomorphism
$$
\{\text{Paths of length $l$ in $Q$ connecting $i,j$}\}
=\Hom_G(X_i, V^{\ttt l}\ttt X_j)
$$
and that $\Hom_G(X_i, V_l \ttt X_j)$ can be described as the space of
``essential paths'' in $Q$, which is naturally a direct summand in
the space of all paths (see \cite{coq-garcia}). The algebra of essential
paths is also known as the preprojective algebra of $Q$.

By McKay correspondence, $Q$ must be an affine Dynkin diagram. 
We denote by $\Delta(Q)$ the  corresponding affine root system; it has a
basis of simple roots $\al_i, i\in I$. We denote 
$$
L(Q)=\bigoplus_{i\in I} \Z\al_i =\Z^I
$$
the corresponding root lattice. It has a natural  bilinear
form given by $(\al_i,\al_i)=2$ and $(\al_i,\al_j)=-n(i,j)$; as is
well-known, this form is positive semidefinite. The kernel of this
form is $\Z\de$, where $\de$ is the imaginary root of $\Delta(Q)$. 

We denote by $s_i\colon L(Q)\to L(Q)$ the reflection around root
$\al_i$. 

Finally, let $K(G)$ be the Grothendieck group of the category $\Rep G$.
It is  freely generated by classes $[X_i], i\in I$; thus, we have a
natural isomorphism 
\begin{equation}\label{e:K(G)}
\begin{aligned}
K(G)&\isoto L(Q)\\
[X_i]&\mapsto \al_i
\end{aligned}
\end{equation}

\section{Quiver $Q_h$}

We will consider a special class of orientations of $Q$. 

\begin{definition}\label{d:height_function}
A {\em height function} $h$ is a map $I\to \Z$ such that 
$h(i)-h(j)=\pm 1$ if $i,j$ are connected by an edge in $Q$, and
$h(i)\equiv  p(i) \pmod 2$. 
\end{definition}

Every height function gives rise to orientation of edges of $Q$: if
$h(j)=h(i)-1$ then all edges connecting $i$ and $j$ are directed
towards $j$:
$$
i\longrightarrow j \quad\text{ if } h(j)=h(i)-1
$$

We will denote by $Q_h$ the quiver given by this orientation. We will
write $i\to j$ if there exists an edge whose tail is $i$ and head is
$j$. The notation $\sum_{j\colon j\to i}$ will mean the sum over all
vertices  $j$  connected with $i$ by an edge $j\to i$; if there are
multiple edges, the corresponding vertex $j$ will be taken more than
once. 

Orientations obtained in this way will be called {\em
admissible} (note: our use of this word is slightly different from the
use in other sources). It is easy to check that if $Q$ has no loops,
then any orientation of $Q$ is admissible. For type $A$, an
orientation is admissible if the total number of clockwise arrows is
equal to the number of counterclockwise ones (which again rules out
type $\widehat{A}_n$, $n$ even, corresponding to $G=\Z_{n+1}$).  It is
also obvious that $Q_h$ has no oriented loops, and that
adding a constant to $h$ gives the same orientation. 

\begin{figure}[ht]
\fig{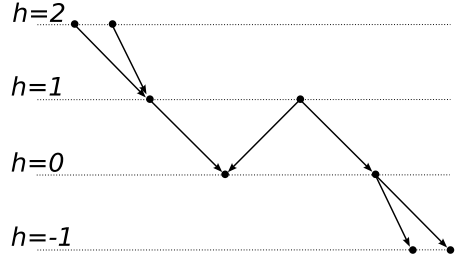}
\caption{Example of  a height function for Dynkin diagram
$\widehat{D}_7$}\label{f:Qh}
\end{figure}

Given a height function $h$, we will draw $Q$ in the plane so that
$h$ is the $y$--coordinate; then all edges of $Q_h$ are directed
down. \firef{f:Qh} shows an example of a height function for
quiver of type $D$.

\begin{definition}\label{d:reflection_functors1}
Let $h$ be a height function on $Q$, and $i\in I$ be a sink in $Q_h$:
there are no edges of the form $i\to j$ (in terms of $h$, it is
equivalent to saying that $h$ has a local minimum at $i$). We define
new height function 
$$
s_i^+h (j)=\begin{cases}
              h(j)+2, & j=i\\
              h(j),   & j\ne\i
           \end{cases}
$$
Similarly, if $i$ is a source, i.e., there are no vertices of the
form $j\to i$ (equivalently, $h$ has a local maximum at $i$), then 
we define 
$$
s_i^-h (j)=\begin{cases}
              h(j)-2, & j=i\\
              h(j),   & j\ne\i
           \end{cases}
$$
\end{definition}
One easily sees that $s_i^+h$, $s_i^-h$ are again height functions;
the quiver $Q_{s_i^\pm h}$ is obtained from $Q_h$ by reversing
orientation of all edges adjacent to $i$. We will refer to $s_i^\pm$
as (elementary) orientation reversal operations. 

The following lemma is known; however, for the benefit of the reader
we included the proof. 

\begin{lemma}\label{l:reflection_functors1}
Any two height functions $h, h'$ can be obtained one from another by
a sequence of orientation reversal operations $s_i^\pm$. 
\end{lemma}
\begin{proof}
Define the ``distance'' between two height functions by \\
$d(h,h')=\sum_I |h(i)-h'(i)|\in \Z_+$. We will show that if $d(h,h')>0$, then
one can
find $s_i^\pm$ such that $s_i^\pm$ can be applied to $h$, and
$d(s_i^\pm h, h')<d(h, h')$; from this the theorem clearly follows. 

Let $I_+=\{i\st h(i)>h'(i)\}$. One easily sees that if $i\in I_+$, and
$j\to i$ in $Q_h$, then $j\in I_+$. Thus, if $I_+$ is non-empty, it 
must contain at least one source $i$ for $Q_h$.  But then $d(s_i^- h,
h')=d(h, h')-2$. 

Similarly, let $I_-=\{i\st h(i)<h'(i)\}$. If $I_-$ is non-empty,
similar
argument shows it  must contain at least one sink $i$ for $Q_h$. 
But then $d(s_i^+ h,h')=d(h, h')-2$. 
\end{proof}

\section{Representations of quivers}
For readers convenience, we recall here the known results about the
relation between representations of quiver $Q_h$ and root system
$\Delta(Q)$ (see, e.g., \cite{drab-ringel}, \cite{kraft-riedtmann}). 

Let $h$ be a height function on $Q$ and $Q_h$ the corresponding quiver.
Consider the category $\Rep Q_h$ of finite-dimensional representations of
$Q_h$; similarly, let $\D^b(Q_h)$ be the bounded derived category of $\Rep
Q_h$. We denote by $K(Q_h)$ the Grothendieck group of $\Rep Q_h$ which can also
be described as the Grothendieck group of $\D^b(Q_h)$. It is well-known that
the category $\Rep Q_h$ is a hereditary abelian  category: 
$$
\Ext^i(X,Y)=0\text{ for all }i>1. 
$$

 For a representation  $M=(M_i)_{i\in I}$ of $Q_h$, we define its
dimension $\dim X\in L(Q)$ by  $\dim M=\sum (\dim M_i)\al_i$ (recall that
$L(Q)=\Z[I]$ is the root lattice of the root system $\Delta(Q)$, see
\seref{s:basic}). 
 
The following theorem summarizes some of the known results about
representations of quivers and root system $\Delta(Q)$.

\begin{theorem}\label{t:reps_of_quivers}
\par\indent
\begin{enumerate}
\item The map $[X]\mapsto \dim X$ gives an isomorphism $K(Q_h)\isoto
L(Q)$. Under this  isomorphism, the bilinear form on $L(Q)$ is identified
with the following   bilinear form on $K(Q_h)$: 
  $$
    (x,y)=\<x,y\>+\<y,x\>
  $$
  where by definition 
  $$
   \<[X],[Y]\>= \dim \RHom (X,Y)=\dim \Hom (X,Y) -\dim \Ext^1(X,Y).
  $$
  \item The set of dimensions of indecomposable modules is exactly the
    set $\Delta_+(Q)$ of positive roots in $\Delta(Q)$. For real roots
    $\al$, there is exactly one  up to isomorphism indecomposable
    module $M_\al$ of dimension $\al$; for
    imaginary root $\al$, there are infinitely many pairwise
    non-isomorphic modules of dimension $\al$. 
  
\end{enumerate}
\end{theorem}

There is also an explicit description of indecomposables in
$\D^b(Q_h)$ (see \cite[Lemma I.5.2]{happel}). 
\begin{theorem}\label{t:indecompos_derived}
  Indecomposable objects in $\D^b(Q_h)$ are of the form $M[n]$, where
  $M$ is an indecomposable object in $\Rep Q_h$, $n\in \Z$. 
\end{theorem}

We will also need reflection functors of Bernstein--Gelfand--Ponomarev.
Recall that if $i$ is a sink in $Q_h$, then one has a natural functor 
\begin{equation}\label{e:reflection_2}
  S_i^+\colon \Rep(Q_h)\to \Rep(Q_{s_i^+h})
\end{equation}
similarly, if $i$ is a source, one has a natural functor 
\begin{equation}\label{e:reflection_3}
  S_i^-\colon \Rep(Q_h)\to \Rep(Q_{s_i^-h})
\end{equation}
(see definition in \cite{drab-ringel},
\cite{kraft-riedtmann}). 

The following result is known and easy to prove, so we skip the
proofs. 
\begin{theorem}\label{t:coxeter2}
\par\indent 
  \begin{enumerate}
  \item Functor $S_i+$ is left exact and $S_i^-$ is right exact.

  We will denote by $RS_i^+, LS_i^-\colon \D^b(Q_h)\to
    \D^b(Q_{s_i^\pm h})$ the  corresponding derived functors. 
  \item
    The functors $RS_i^+$, $LS_i^-$ are equivalences of categories
    $\D^b(Q_h)\to\D^b(Q_{s_i^\pm h})$, which induce the usual
    reflections $s_i$ on the  Grothendieck group: 
    $$
    \dim RS_i^+ (X)=s_i(\dim X), \quad \dim LS_i^- (X)=s_i(\dim X)
    $$
    \item If $i,j$ are not neighbors in $Q$, then $RS_i^+,
    RS_j^+$ commute (i.e., compositions in different orders are
    isomorphic) and similarly for $LS_i^-$. 
  \end{enumerate}
\end{theorem}

In particular, for a given height function $h$ let $s^+_{i_1}\dots
s^+_{i_r}$ be a sequence of elementary orientation
reversal operations such that  
$$
s_{i_1}^+\dots s_{i_r}^+ (h) =h+2
$$
One easily sees that this condition is equivalent to requiring that every
index $i\in I$ appears in the sequence $\{i_1,\dots, i_r\}$ exactly once;
it follows from \leref{l:reflection_functors1}  that for every height
function $h$, such sequences of orientation reversal operations exist. 
For such a sequence,  the corresponding element of the  Weyl group
\begin{equation}\label{e:coxeter}
c^+_h=s_{i_1}\dots s_{i_r}
\end{equation}
is called the {\em Coxeter element}, and
the corresponding composition of reflection functors 
\begin{equation}\label{e:C+}
RC^+_h=RS_{i_1}^+\dots RS_{i_r}^+\colon \D^b(Q_h)\to\D^b(Q_{h+2})
\end{equation}
will be called the {\em Coxeter functor}. Note that since $Q_h\simeq
Q_{h+2}$ as a quiver, we can consider$RC_h^+$ as an autoequivalence of
$\D^b(Q_h)$. 

It is easy to show (see \cite{drab-ringel}, \cite{shi}) that the
Coxeter element $c^+_h$ only depends on $h$ and not on the choice of
the sequence $i_1,\dots,i_r$; moreover, the proof of this only uses
the fact that $s_i, s_j$ commute if $i,j$ are not connected in $Q$
and does not use the braid relations. Thus, by \thref{t:coxeter2},
this implies that up to an isomorphism, $RC^+_h$ also depends only on
$h$; this justifies the notation $RC_h^+$.

Similarly, we can define  functors 
\begin{equation}\label{e:coxeter2}
LC^-_h\colon \D^b(Q_h)\to\D^b(Q_{h-2})\simeq\D^b(Q_h)
\end{equation}
using sequences of orientation reversals  $s_{i_1}^-\dots s_{i_r}^- h
=h-2$; the corresponding element of the Weyl group will be denoted by
$c_h^-$. As before, it can be shown that $LC_h^-$, $c_h^-$ only depend
on $h$.

For readers familiar with the theory of Auslander--Reiten sequences
(see \cite{auslander-reiten}, \cite{happel}), we add
that the category $\Rep Q_h$ has Auslander--Reiten sequences, and the
Auslander--Reiten functor $\tau$ is given by $\tau=C_h^-$.

\section{Equivariant sheaves}
In this section, we introduce the main object of this paper, the
category of equivariant sheaves on $\PP$. Most of  the results of
this section  are well-known and given here only for the convenience
of references. \leref{l:frobenius} does not seem to be easily
available in the literature, but is very easy to prove. 

Let $V^*$ be the dual of the tautological representation $V$ of
$\SU$. Since $G$ is a finite subgroup in $\SU$, it naturally acts on
$\PP=\mathbb{P}(V^*)$, and the structure sheaf $\OO$ has a standard
$\SU$- (and thus $G$-) equivariant structure. Moreover, all twisted
sheaves $\OO(n)$ also have a standard equivariant structure, so that
the space of global sections $\Ga(\OO(n))$ is a representation of 
$G$: 
\begin{equation}\label{e:H0}
\Ga(\OO(n))=\begin{cases}
                  S^nV=V_n, &n\ge 0\\
                  0, &n<0
                  \end{cases}
\end{equation}
Similarly, the higher homology spaces are naturally representations of $G$: it
is well-known that $H^i(\PP, \OO(n))=0$ for $i>1$, and 
\begin{equation}\label{e:H1}
H^1(\PP,\OO(n))=\begin{cases}
                  S^{-n-2}V^*=V^*_{-n-2}, &n\le -2\\
                  0, &n\ge -1
                  \end{cases}
\end{equation}

Let $\Gqcoh$, $\Gcoh$ be the categories of $G$--equivariant
quasi-coherent (respectively, coherent) $\OO$--modules on $\PP$
(see, e.g., \cite[Section 4]{bkr} for definitions). Note that we are
considering isomorphisms $\la_g\colon \F\to g^*\F$ as part of the
structure of the $G$--equivariant sheaf. For brevity, we will denote 
morphisms and $\Ext$ groups in $\Gqcoh$ by $\GHom(X,Y)$, $\GExt(X,Y)$.
For an equivariant  sheaf $\F$ we will denote 
$$
\F(n)=\OO(n)\ttt_\OO \F
$$
with the obvious equivariant structure. Similarly, for a
finite-dimensional representation $X$ of $G$, we denote 
$$
X(n)=\OO(n)\ttt_\CC X.
$$

We list here some of the basic properties of equivariant
sheaves; proofs can be found in \cite[Section 4]{bkr}.
\begin{theorem}\label{t:Gcoh}
\par\noindent
  \begin{enumerate}
    \item $\Gqcoh$ is an abelian category, and a sequence $0\to
    \F_1\to \F_2\to \F_3\to 0$ is exact in $\Gqcoh$ iff it is exact in
    $\Qcoh$. 

  \item For any $\F,\G\in \Gqcoh$, the space $\Hom_\OO(\F,\G)$ is a
    representation of $G$, and $\GHom(\F,\G)=(\Hom_\OO(\F,\G))^G$.
    Similarly, \\
    $\GExt^i(\F,\G)=(\Ext^i_\OO(\F,\G))^G$; in
    particular, $\GExt^i(\F, \G)=0$ for $i>1$.
    
   \item For any $\F,\G\in \Gcoh$, the spaces  $\GHom(\F,\G)$,
    $\GExt^1(\F,\G)$ are finite-dimensional. 

  \item For any $\F,\G\in \Gcoh$, one has 
  \begin{align*}
  \GHom(\F, \G(n))=\GHom(\F(-n),\G)=0\quad \text{for }n\ll 0,\\
  \GExt^1(\F, \G(n))=\GHom(\F(-n),\G)=0\quad \text{for }n\gg 0,
  \end{align*}
 
\end{enumerate}
\end{theorem}

As an immediate corollary of this, we see that the category $\Gcoh$ has the 
Krull--Schmidt property: every object $\F\in \Gcoh$ can be written as a direct
sum of indecomposable modules, and the multiplicities do not depend on the
choice of such decomposition. It is also a hereditary category (recall that a
category is called hereditary if $\Ext^2(A,B)=0$ for any objects $A,B$; in
particular, this implies that a quotient of an injective object is injective,
as can be easily seen from the long exact sequence of $\Ext$ groups).  

We say that $\F\in\Gcoh$ is locally free if it is locally free as an
$\OO$-module. 
\begin{theorem}\label{t:locally_free}
\par\noindent
\begin{enumerate}
  \item Every locally free $\F\in\Gcoh$ is isomorphic to direct sum of
    sheaves of the  form 
    \begin{equation}
    X_i(n)=\OO(n)\ttt X_i, \quad X_i \text{ -- irreducible representation of
$G$}    \end{equation}
  \item If $\F\in\Gcoh$ is locally free, then the functor $\ttt \F\colon
    \Gcoh\to \Gcoh$ is exact. 
   
 \item Every coherent $G$-equivariant sheaf has a resolution consisting of
    locally free sheaves. 
 \item \textup{(}Serre Duality\textup{)} 
   For any two locally free sheaves $\F, \G$ we have isomorphisms 
   $$
   \GExt^1(\F, \G(-2))=\Ext^1(\F(2),\G)=\GHom(\G, \F)^* 
   $$
\end{enumerate}
\end{theorem}

It immediately follows from computation of homology of $\OO(n)$ as coherent
sheaves that 
\begin{equation}\label{e:hom1}
  \GHom ( X_i(n), X_j(k))
  =\begin{cases}
     \Hom_G( X_i, V_{k-n}\ttt X_j), &k\ge n\\
     0, & k<n
      \end{cases}
\end{equation}
\begin{equation}\label{e:hom2}
  \GExt^1 (X_i(n), X_j(k))
  =\begin{cases}
     \Hom_G(X_i, V^*_{n-k-2}\ttt X_j), &k\le n-2\\
     0, & k\ge n-1
      \end{cases}
\end{equation}

Finally, we note that the action of $(-I)\in \SU$ gives a decomposition
of $\Gcoh$ into even and odd part. Namely, since $-I$ acts trivially on
$\PP$, structure of $G$-equivariant sheaf gives an isomorphism
$\la_{-I}\colon  \F\to (-I)^*\F=\F$. We define 
$$
\Gcoh_p=\{\F\in\Gcoh\st \la_{-I}=(-1)^p\}, 
\quad p\in \Z_2
$$

From now on, the main object of our study will be the category 
\begin{equation}\label{e:C}
\C=\Gcoh_0.
\end{equation}
In particular, a locally free sheaf $X_i(n)\in 
\C$ iff $n+p(X_i)\equiv 0 \mod 2$, where $p(X_i)$ is the parity defined by
\eqref{e:parity}. 

One easily sees that $\C$ is a full subcategory in $\Gcoh$ closed under
extensions. Equivalently, it can be described as follows. 

\begin{lemma}
Let $\Gbar=G/\{\pm I\}$.Then $\C$ is naturally equivalent to the
category $\Gbarcoh$ of $\Gbar$-equivariant coherent sheaves on $\PP$.
\end{lemma}

We will denote the bounded  derived category of $\C$ by  $\D^b(\C)$.

For future use, we will also need an analogue of induction functor.
Recall that for any finite subgroup $H\subset G$, we have an
induction functor $\Ind_H^G\colon \Rep H\to \Rep G$. In particular,
for  $H=\{1\}$ and the trivial representation $\CC$ it gives the
regular representation of $G$:
\begin{equation}\label{e:regular_rep}
R=\Ind_{\{1\} }^G\CC\simeq \bigoplus_{i\in I}d_i X_i,
\qquad d_i=\dim X_i
\end{equation}
As any representation of $G$, $R$ can be decomposed into even and odd
part (cf. \eqref{e:even-odd}): $R=R_0\oplus R_1$, where 
\begin{equation}\label{e:R_p}
R_p=\bigoplus_{i\in I_p}d_iX_i
\end{equation}

It is easy to check that $R_0$ is exactly the regular representation
of $\Gbar$. 

Similarly, we have an induction functor from the category of
$H$-equivariant sheaves to $G$--equivariant sheaves; in particular,
for $H=\{1\}$, we get a functor $\Coh\to\Gcoh$. However, it will be
more convenient to consider the functor 
\begin{equation}\label{e:ind}
  \begin{aligned}
  \Ind\colon \Coh&\to\C=\Gbarcoh\\
            \F&\mapsto \bigoplus_{g\in\Gbar}g^*\F 
  \end{aligned}
\end{equation}
The following lemma lists some of the properties of this functor; the
proof is straightforward and left to the reader.
\begin{lemma}\label{l:frobenius}
Let $\Ind\colon \Coh\to\C$ be defined by \eqref{e:ind}. Then for even
$n$, $\Ind \OO(n)$ is naturally isomorphic as an 
    equivariant sheaf to $R_0\ttt \OO(n)$; for odd $n$, $\Ind \OO(n)$
      is    naturally isomorphic as an 
    equivariant sheaf to $R_1\ttt \OO(n)$.
  
\end{lemma}

\section{Auslander--Reiten relations}
The following result will play a key role in the study of the category
$\Gcoh$. 

\begin{theorem}\label{t:auslander_reiten}
For any $n\in \Z$, $i\in I$, there is short exact sequence in $\Gcoh$: 
\begin{equation}\label{e:aus-reiten1}
0\to X_i(n)\to \sum_j X_j(n+1)\to X_i(n+2)\to 0
\end{equation}
where the sum is over all neighbors $j$ of $i$ in $Q$. 
\end{theorem}
We will call sequences of the form \eqref{e:aus-reiten1}  {\em
Auslander--Reiten sequences}.

\begin{proof}
Notice first that we have a short exact sequence of $\OO$-modules:
$$
0\to \OO\to \OO(1)\ttt V\to \OO(2)\to 0
$$ 
Let us tensor it with $X_i(n)$. Since tensoring with a locally free sheaf
is an exact functor (\thref{t:locally_free}), this  gives
short exact sequence 
\begin{equation}\label{e:aus-reiten2}
0\to X_i(n)\to  (V\ttt X_i)(n+1)\to X_i(n+2)\to 0
\end{equation}
Since $V\ttt X_i=\sum_{j} X_j$, we get the statement of the theorem. 
\end{proof}

\begin{remark}
For readers familiar with the general theory of Auslander--Reiten (AR)
sequences (see \cite{auslander-reiten} for the overview of the
theory), we point out that our use of this name is justified: it can
be shown, using Serre duality that these sequences do satisfy the
usual definition of AR sequences (see \cite{reiten-bergh}). This, in
particular, shows that the Auslander--Reiten functor $\tau$ for locally
free equivariant sheaves is given by 
$$
\tau(\F)=\F(-2).
$$
However, we are not going to use the general theory of AR sequences in
this paper.
\end{remark}
\begin{corollary}
  For any $n\in \Z, i\in I$ we have the following relations in
  the Grothendieck group $K(\Gcoh)$
    \begin{equation}\label{e:aus-reiten3}
    [X_i(n)]-\sum_{j}[X_j(n+1)]+[X_i(n+2)]=0
    \end{equation}
  where the sum is over all neighbors $j$ of $i$ in $Q$.
\end{corollary}
Later we will study the Grothendieck group $K(\Gcoh)$ in more detail,
in particular showing that relations \eqref{e:aus-reiten3} is the full
set of relations among classes $[X_i(n)]$ (see
\coref{c:generators_of_K(C)}).

\section{Indecomposable objects in $\Gcoh$}
In this section, we describe the indecomposable sheaves in $\Gcoh$.
Again, these results are not new; they  follow from results of
\cite{lenzing}, \cite{schiffmann} who considers more general setting of
``non-commutative projective curves''. So this section is included just for
the reader's convenience.

We start by describing torsion sheaves. 

For every point $x\in \PP$, $n\in \Z_+$, denote $\OO_{n[x]}=\OO/m_x^n$,
where $m_x$ is the ideal of functions vanishing at $x$. This sheaf is
supported at point $x$; choosing a local coordinate $z$ at $x$, we can
identify  the stalk of this sheaf at $x$ with $\CC[z]/z^n\CC[z]$. It
is well known that sheaves $\OO_{n[x]}$ are indecomposable, and that
every  coherent $\OO$-module is isomorphic to a direct
sum of a locally free sheaf and sheaves of the form $\OO_{n[x]}$. 

Assume now that  $x\in \PP$ is such that $\Stab_G(x)=\{\pm I\}$ and
thus $\Stab_\Gbar x=\{1\}$ (recall that $\Gbar=G/\{\pm I\}$). Such
points will be called  {\em generic}. Define
\begin{equation}\label{e:torsion_sheaf}
\OO_{n[Gx]}=\Ind \OO_{n[x]}=\bigoplus_{g\in
\Gbar}\OO_{n[gx]}
\end{equation}
where the functor  $\Ind$ is defined by \eqref{e:ind}.
This sheaf is supported on the orbit $Gx$ and has a canonical
$\Gbar$--equivariant structure which can be lifted to a
$G$-equivariant structure. It is easy to see that the space of global
sections $\Ga(\OO_{n[Gx]})$ is isomorphic to $n$ copies of the regular
representation of $\Gbar$.

This construction can be generalized to non-generic points.

\begin{theorem}\label{t:torsion}
Let $x\in \PP$ be non-generic: $H=\Stab_\Gbar x\ne \{1\}$. Then 

\begin{enumerate}
\item $H$ is  a cyclic group. 
\item For suitable choice of the generator $\si$ of $H$, its action in
  the one-dimensional  vector space $\OO_{[x]}=\OO/m_x$ is given by
  $\si\mapsto e^{2\pi i /|H|}$.  
\item Consider the completed local ring  $\widehat{\OO}_{\PP,x}$;
  choosing a local coordinate $z$ at $x$ identifies this algebra with
  $\CC[[z]]$. It has a natural action of $H$. Then we have a natural
  equivalence of categories
  \begin{align*}
  &\{\Gbar\text{-equivariant coherent sheaves supported on }Gx\}\\
  &\quad \simeq
   \{H\text{-equivariant finite-dimensional modules over }
     \widehat{\OO}_{\PP,x}\}             
  \end{align*}
  given by 
  $$
  \F\mapsto\text{stalk of $\F$ at $x$}
  $$

\end{enumerate}
\end{theorem}
The proof of this theorem is straightforward and  left to the reader.

\begin{definition}
Let $H=\Stab_\Gbar x$, and $Y$ --- a finite-dimensional representation
of $H$. Then we denote by $\OO_{n[Gx],Y}$ the $\Gbar$-equivariant sheaf
supported on the orbit  $Gx$ and whose
    stalk at $x$ is isomorphic as an $H$-module to   $Y\otimes
    (\OO_{\PP, x}/m_x^n)$. 

\end{definition}
It is easy to show that the space of global sections of $\OO_{n[Gx],Y}$
is isomorphic to 
  $$
  \Ind_H^G (Y\otimes (\OO/m_x^n)).
  $$

\begin{theorem}\label{t:indecomposable}
\par\noindent
\begin{enumerate}
  \item
    The following is the full list of indecomposable objects in
  $\C$: 
    \begin{itemize}
      \item Torsion sheaves $\OO_{n[Gx]}$, $x\in \PP$,
      $\Stab_{\Gbar}x=\{1\}$, $n\in \Z_+$
      \item Torsion sheaves $\OO_{n[Gx],Y}$, $x\in \PP$,
        $H=\Stab_{\Gbar }x\ne\{1\}$, 
        $Y$-- an irreducible representation of $H$.
      \item Free sheaves $X_i(n)$, $n\in \Z$, $i\in I$, $n+p(i)\equiv
      0\mod $.  
    \end{itemize}
  \item Indecomposable objects in $\D^b(\C)$ are of the form $M[n]$,
    $M$---an indecomposable object in $\C$, $n\in \Z$. 
  \item There are no injective and no projective indecomposable objects in
    $\C$. 
\end{enumerate}
\end{theorem}
\begin{remark}
  Of course, sheaves $\OO_{n[Gx]}$, $H=\Stab_{\Gbar}x=\{1\}$, can be
  considered as special case of sheaves $\OO_{n[Gx],Y}$. However, we
  choose to list cases $H=\{1\}$ and $H\ne\{1\}$ separately. 
\end{remark}

\begin{proof}
It is known that every coherent sheaf $\F$ has a maximal torsion
subsheaf $\TT$ so that $\F/\TT$ is locally free, and short exact
sequence $0\to\TT\to \F\to \F/\TT\to 0$ splits. If we additionally
assume that $\F$ is a $\Gbar$-equivariant sheaf, then $\TT$ (and
thus $\F/\TT$) must also be $\Gbar$-equivariant. It easily follows from
\thref{t:Gcoh} that then $0\to\TT\to\F\to\F/\TT\to 0$ splits in the
category of $G$-equivariant sheaves; thus, every coherent
$G$-equivariant sheaf is a direct sum of a free sheaf and a torsion
sheaf. 

By \thref{t:locally_free}, indecomposable free sheaves are of the form
$X_i(n)$. For torsion sheaves, it is easy to see that an indecomposable
torison sheaf must be supported on an orbit $Gx$. Now the
classification of indecomposable torsion sheaves follows from
\thref{t:torsion}.

  Classification of indecomposable objects in $\D^b$ follows from 
\cite[Lemma I.5.2]{happel}. His argument was given for the category of
representations of a hereditary algebra, but actually works in any
hereditary category:  it only  uses that a quotient of an injective
object is injective, which easily follows from long exact sequence of
$\Ext$ spaces and the fact that $\Ext^2(\F,\G)=0$ for any $\F,\G$.
\end{proof}

\section{Auslander--Reiten quiver $\Qhat$}
We can now define a combinatorial structure which will play a crucial
role in our paper. Let $\Qhat$ be the set of isomorphism classes of
locally free indecomposable objects in $\C$;
by \thref{t:indecomposable}, we can write 
\begin{equation}\label{e:qhat}
\Qhat=\{(i,n)\st n\in \Z, i\in I, n+p(i)\equiv 0 \mod 2\}
\end{equation}
We will frequently use the notation 
$$
X_q=X_i(n)=\OO(n)\ttt X_i, \quad q=(i,n)\in \Qhat.
$$

We turn $\Qhat$ into a  quiver by defining, for $q=(i,n)$,
$q'=(j,n+1)$, 
\begin{align*}
&(\text{Number of edges }q\to q')=\dim \Hom_\C(X_q,X_{q'})\\
&\qquad =\dim \Hom_G(X_i, X_j\ttt V)=n(i,j)
\end{align*}
(recall that $n(i,j)$ is the number of edges between $i$ and $j$ in
$Q$ and $V=\Ga(\PP,\OO(1))$ is the tautological representation of
$\SU$).

Note that the edges only connect vertices with $n$ differing by one,
and all edges are directed ``up'' (i.e., so that $n$ increases). The
figure below shows an example of the quiver $\Qhat$ when $G$ is of
type $\widehat{D_7}$. 
\begin{figure}[ht]
\fig{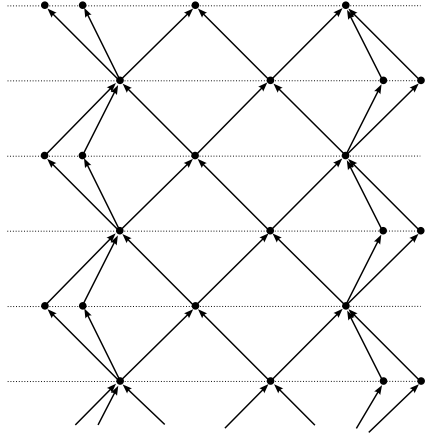}
\caption{Auslander--Reiten quiver $\Qhat$ for $Q=\widehat{D_7}$}
\end{figure}

\begin{remark}
Again, for users familiar with the Auslander--Reiten theory, we note that
$\Qhat$ is exactly the stable AR quiver of  the category $\Gcoh$. As
was mentioned in the introduction, it is the picture of this quiver in
Ocneanu's work (who defined it purely combinatorially, without
reference to the category of sheaves) that inspired the current paper. 
\end{remark}

The following easy lemma is left as an exercise to the reader. 
\begin{lemma}\label{l:qhat_connected}
$\Qhat$ is a connected quiver. 
\end{lemma}
This is actually the reason for considering the even part $\C=\Gcoh_0$
rather than all of $\Gcoh$: otherwise, the resulting quiver would not be
connected.

Note that  we have a natural pairing 
$$
\Hom_G(X_i, X_j\ttt V)\ttt\Hom_G(X_j, X_i\ttt V)\to \CC
$$
given by composition
$$
X_i\xxto{\ph}X_j\ttt V\xxto{\psi\ttt 1}
 X_i\ttt V^{\ttt 2}\xxto{1\ttt ev}X_i 
$$
where $ev\colon V^{\ttt 2}\to \CC$ is a $G$-equivariant pairing (which
is defined uniquely up to a constant).  This
gives rise to duality 
\begin{equation}\label{e:pairing}
\Hom_\C(X_i(n), X_j(n+1))
=\Hom_\C(X_j(n+1), X_i(n+2))^*
\end{equation}

As usual in quiver theory, we define path algebra of $\Qhat$ by 
\begin{equation}\label{e:path_algebra}
\begin{aligned}
&\Path(q,q')=\text{Span(paths from $q$ to $q'$ in $\Qhat$)},\\
&\Path=\bigoplus_{q,q'\in \Qhat}\Path(q,q')
\end{aligned}
\end{equation}
One easily sees that $\Path( (i,k), (j,l))$ coincides with the space of
paths in $Q$ of length $l-k$ between $j$ and $j$ (if $l-k\ge 0$; otherwise,
$\Path((i,k), (j,l))=0$ ).

We also define an analogue of preprojective algebra by 
\begin{equation}\label{e:preproj_algebra}
\begin{aligned}
&\Pi_\Qhat=\Path/(\theta_q)\\
&\theta_q=\sum_{e\colon q\to q'}  e_i e^i\colon q\to q(2)
\end{aligned}
\end{equation}
where for $q=(i,n)$, we denote $q(2)=(i,n+2)$, the sum is over all edges
originating at $q$, and $e_i, e^i$ are dual bases in $\Hom_\C(X_q, X_{q'})$
and $\Hom_\C(X_{q'}, X_{q(2)})$ with respect to the pairing
\eqref{e:pairing}. (This algebra is isomorphic to
Ocneanu's algebra of essential paths, or, in the terminology of AR
quivers,  to  the mesh algebra of $\Qhat$ considered as a polarized
translation quiver.)

\begin{theorem}\label{t:homs=paths}
Let $q,q'\in \Qhat$. Then 
$$
\Hom_\C(X_q,X_{q'})=\Pi_\Qhat(q,q')
$$ 
\end{theorem}
\begin{proof}
Let $\Pi_Q$ be the preprojective algebra of the graph $Q$, and $\Pi_Q^l$
subspace of paths of length $l$. Then for $q=(i,k)$, $q'=(j,l)$ it follows
from the definition that for $l<k$ we have 
$$
\Pi_\Qhat (q,q')=0=\Hom_\C(X_i(k),X_j(l))
$$ 
and for $l\ge k$, 
$$
\Pi_\Qhat (q,q')=\Pi_Q^{l-k}(X_i,X_j)=
\Hom_G(X_i,V^{\ttt (l-k)}\ttt X_j)
=\Hom_\C(X_i(k),X_j(l)).
$$
\end{proof}

This theorem shows that morphisms between free indecomposable
modules in $\C$ are described by essential paths in $\Qhat$, so the
structure of the subcategory $Free(\C)$, consisting of locally free sheaves,
can be recovered from $\Qhat$. Note however that  realization of
$Free(\C)$ in terms of equivariant sheaves gives much more than
$\Hom$ spaces: it embeds $Free(\C)$ into an
abelian category with sufficiently many injectives, allowing one to
define $\Ext$ functors and derived category $\D^b(\C)$. Both of
these would be difficult to define in terms of $\Qhat$.

\section{Equivalence of categories}
Let $Q$, $\Qhat$ be as defined above, and let $h$ be a height
function on $Q$ (see \deref{d:height_function}). Recall that such 
a height function defines a choice of orientation on $Q$; the
corresponding quiver is denoted $Q_h$. Then this height function
defines an embedding 
\begin{equation}\label{e:embedding}
\begin{aligned}
i_h\colon Q_h^{opp}& \to \Qhat\\
    i&\mapsto (i,h_i)
\end{aligned}
\end{equation}
where $Q_h^{opp}$ is the quiver obtained from $Q_h$ by reversing all
arrows. The figure below shows an example of such an embedding.

\begin{figure}[ht]
\fig{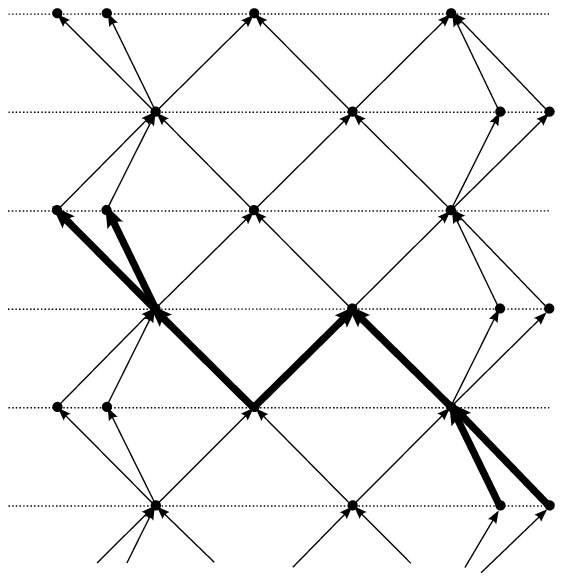}
\end{figure}

\begin{remark}
Note that once we have embedded $Q_h\subset \Qhat$, we can 
identify the quiver $\Qhat$ with $\Z Q_h$ (see, e.g., \cite{ringel1}
for the definition).
\end{remark}

From now on, we will consider $Q_h^{opp}$ as a subset of $\Qhat$,
omitting $i_h$ in our notation. 

%

\begin{definition}
Let $q_1,q_2\in \Qhat$.  We say that $q_1$ is a predecessor of $q_2$
(notation: $q_1\prec q_2$) if there exists a path $q_1\to
\bullet\to\dots\to \bullet \to q_2$ in $\Qhat$. In this case, we also
say that $q_2$ is a successor of $q_1$ and write $q_2\succ q_1$. 
\end{definition}

\begin{lemma}\label{l:technical1}
  Let $q=(i,n)\in \Qhat$. Then 
  \begin{enumerate}
    \item $n>h_i \iff (q\notin Q_h, q\succ q'$ for some $q'\in Q_h)$. In
    this case,  we will say that $q$ is above $Q_h$ and write
    $q\succ Q_h$.
    \item $n<h_i \iff (q\notin Q_h, q\prec q'$ for some $q'\in Q_h)$. In
      this case,  we will say that $q$ is below $Q_h$ and write
      $q\prec Q_h$.
  \end{enumerate}
\end{lemma}
\begin{lemma}\label{l:vanishing_of_homs}
\par\indent
  \begin{enumerate}
    \item If $q_1\succeq Q_h,  q_2\prec Q_h$, then
        $\Hom_\C(X_{q_1},X_{q_2})=0$
    \item If   $q_1\preceq Q_h,  q_2\succeq Q_h$, then
        $\Ext^1_\C(X_{q_1},X_{q_2})=0$. In particular, if $q_1, q_2\in
        Q_h$, then $\Ext^1_\C(X_{q_1},X_{q_2})=0$. 
  \end{enumerate}
\end{lemma}
\begin{proof}
  It follows from $\Hom_\C(X_{q_1},X_{q_2})=\Pi_{\Qhat}(X_{q_1},X_{q_2})$
  that  $\Hom_\C(X_{q_1},X_{q_2})$ can be nonzero only if $q_1\preceq q_2$.
  But  then  $q_2\succeq q_1\succeq Q_h$, which proves part (a). 
  
  Part (b) follows from part (a) and Serre duality
  (\thref{t:locally_free}) 
\end{proof}

\begin{lemma}\label{l:generators}
The sheaves $X_q$, $q\in Q_h$, generate the category $\D^b(\C)$ as a
triangulated category: the smallest full triangulated subcategory in
$\D^b(\C)$ containing all $X_q$, $q\in Q_h$, is $\D^b(\C)$. 
\end{lemma}
\begin{proof}
Let $\D$ be the smallest triangulated subcategory in
$\D^b(\C)$ containing all $X_q$, $q\in Q_h$. If $i$ is a sink for $Q_h$,
then it follows from AR exact sequence \eqref{e:aus-reiten1} that
$X_i(h_i+2)\in \D$. Thus, all $X_q$, $q\in Q_{s_i^+ h}$, are in $\D$.
Similarly, if $i$ is a source for $Q_h$, then all $X_q$, $q\in Q_{s_i^-
h}$, are in $\D$. By \leref{l:reflection_functors1}, this implies that
all $X_q$, $q\in \Qhat$, are in $\D$, so $\D$ contains all locally free
equivariant sheaves. Since every coherent sheaf admits a locally free
resolution, this completes the proof. 
\end{proof}

Define now a functor 
\begin{equation}\label{e:Ph}
\begin{aligned}
\Ph_h\colon \C\to \Rep(Q_h)
\end{aligned}
\end{equation}
by 
$$
\Ph_h(\F)(i)=\Hom_\C(X_i(h_i), \F)
$$
and for an edge $e\colon j\to j$ in $Q_h$, we define the corresponding
map by  
\begin{align*}
\Hom_\C(X_j(h_j), \F)&\to \Hom_\C(X_i(h_i), \F)\\
                  \ph&\mapsto \ph\circ h(e)
\end{align*}
where $h(e)\colon (i,h_i)\to (j,h_j)$ is the edge in $\Qhat$ corresponding to
$e$ under embedding \eqref{e:embedding}.

\begin{lemma}\label{l:Ph}
\par\indent
\begin{enumerate}
\item $\Ph_h$ is left exact.
\item $R^i\Ph_h=0$ for $i>1$, and 
$$
(R^1\Ph_h(\F))(i)=\Ext_\C^1(X_i(h_i),\F)
$$
%
\end{enumerate}
\end{lemma}
\begin{proof}
Follows from left exactness of functor $\Hom_\C(X_i(h_i), -)$ and definitions. 
\end{proof}

\begin{example}\label{x:standard}
Let $h$ be the ``standard'' height function: $h(i)=0$ for $i\in I_0$,
$h(i)=1$ for $i\in I_1$. Then the map $\Ph_h$ can also be described
as follows. Consider the functor $\Rep Q_h\to \Rep G$ defined by
$\{V_i\}\mapsto \bigoplus V_i\ttt X_i$. Then the composition
$$
\C\xxto{\Ph_h} \Rep Q_h\to \Rep G
$$
is given by 
$$
\F\mapsto \Ga(\F)\oplus \Ga(\F(-1)).
$$

Indeed, for each representation $M$ of $G$ we have 
$$
M\simeq\bigoplus \Hom_G(X_i,M)\ttt X_i
$$
Applying it to $\Ga(\F)=\Hom_\OO(\OO,\F)$ we get 
$$
\Ga(\F)=\bigoplus_{i\in I} \Hom_\OO(X_i\ttt \OO,\F)^G\ttt X_i
=\bigoplus_{i\in I} \GHom (X_i\ttt \OO,\F)\ttt X_i
$$
Since $\F\in \C=\Gcoh_0$, we have $\GHom (X_i\ttt \OO,\F)=0$ for
$i\in I_1$; thus, 
$$
\Ga(\F)=\bigoplus_{i\in I_0} \GHom (X_i\ttt \OO,\F)\ttt X_i
$$
Similar argument shows that 
$$
\Ga(\F(-1))=\bigoplus_{i\in I_1} \GHom (X_i\ttt \OO(1),\F)\ttt X_i
$$
Thus, 
$$
\Ga(\F)\oplus \Ga(\F(-1))=
\left(\bigoplus_{i\in I_0} \Hom_\C (X_i,\F)\ttt X_i\right)
\oplus
\left(\bigoplus_{i\in I_1} \Hom_\C (X_i(1),\F)\ttt X_i\right)
$$
as desired. 
\end{example}

Since $\Ph_h$ is left exact, we can define the associated derived functor
$R\Ph_h\colon \D^b(\C)\to \D^b(\Rep Q_h)$. 
The following two theorems are the main results of this paper. 

\begin{theorem}\label{t:main2}
For any height function $h$, the functor $R\Ph_h\colon \D^b(\C)\to
\D^b(Q_h)$ defined by \eqref{e:Ph} is an equivalence of triangulated
categories. 
\end{theorem}
\begin{proof}
It follows from \leref{l:vanishing_of_homs} that the object 
$$
T=\bigoplus_{q\in Q_h}X_q
$$
satisfies $\Ext^1(T,T)=0$. By \leref{l:generators}, direct summands of
 $T$ generate $\D^b(\C)$ as a triangulated category. Therefore, $T$ is
the tilting object in $\D^b(\C)$ as defined in \cite{geigle-lenzing}.
Now the statement of the theorem follows from the general result of
\cite{geigle-lenzing}.
\end{proof}

\begin{theorem}\label{t:main1}
Let $i$ be a sink for $h$, and $S_i^+$ the reflection functor  defined by 
\eqref{e:reflection_2}. Then the following diagram is commutative:
$$
\xymatrixrowsep{4pt}
\xymatrixcolsep{40pt}
\xymatrix{
                  &\D^b(Q_{s_i^+h})\\
\D^b(\C)\ar[ur]^{R\Ph_{s_i^+h}}\ar[dr]^{R\Ph_{h}}&\\
                 &\D^b(Q_{h})\ar[uu]_{RS_i^+}
} 
$$
and similarly for $S_i^-$.  
\end{theorem}
\begin{proof}

Let us first prove that if $q\succ Q_h$, then 
$$
R\Ph_{s_i^+h}(X_q)=RS_i^+\circ R\Ph_h (X_q).
$$ 
Indeed, in this case it follows
from \leref{l:Ph}, \leref{l:vanishing_of_homs} that $R^i\Ph_h(X_q)=0$
for $i>0$,
so $R\Ph_h(X_q)=\Ph_h(X_q)$; similarly, 
$R\Ph_{s_i^+h}(X_q)=\Ph_{s_i^+h}(X_q)$.

Let $i$ be a sink for $Q_h$.  Then we have the Auslander--Reiten exact
sequence
$$
0\to X_i(h_i)\to \bigoplus_{j} X_j(h_j)\to X_i(h_i+2)\to 0
$$
Applying to this sequence functor $\Hom_\C(-, X_q)$ and using
\leref{l:vanishing_of_homs} which gives  vanishing of
$\Ext^1$,  we get a short exact
sequence 
$$
0\to \Hom_\C(X_i(h_i+2), X_q)\to \bigoplus_{j} \Hom_\C(X_j(h_j), X_q)\to
\Hom_\C(X_i(h_i), X_q)\to 0
$$
Comparing this with definition of $S_i^+$, we see that 
$\Ph_{s_i^+h}(X_q)=S_i^+(\Ph_h(X_q))$. 

Thus, we have shown that for $q\succ Q_h$, we
have $R\Ph_{s_i^+h}(X_q)=RS_i^+\circ R\Ph_h (X_q)$.

To complete the proof, we now use the following easy result.

\begin{lemma}\label{l:technical4}
   Let $\D$, $\D'$ be triangulated categories, $\Ph_1,\Ph_2\colon
  \D\to   \D'$ triangle functors, and $\al\colon
  \Ph_1\to\Ph_2$ a morphism of  functors. Assume that  there is a
  collection of objects $X_q\in \D$ which generate $\D$ as a
  triangulated category in the sense of \leref{l:generators} and such
  that for each $X_q$, 
  $$
  \al\colon \Ph_1(X_q)\to \Ph_2(X_q)
  $$
  is an isomorphism. Then $\al$ is an isomorphism of functors. 
\end{lemma}  

Applying this lemma to $\D=\D^b(\C)$, $\D'=\D^b(Q_{s_i^+h})$,
$\Ph_1=R\Ph_{s_i^+ h}$, $\Ph_2=RC_i^+\circ R\Ph_h$ and the 
collection of objects $X_q, q\succ Q_h$ (which
generate $\D^b(\C)$ by \leref{l:generators}), we get the statement of
the theorem. 
\end{proof}

\begin{corollary}\label{c:coxeter}
We have the following commutative diagram
$$
\xymatrixrowsep{20pt}
\xymatrixcolsep{30pt}
\xymatrix{
  \D^b(\C)\ar[r]^-{R\Ph_h}                       
                      &\D^b(Q_h)\simeq \D^b(Q_{h+2})\\
  \D^b(\C)\ar[r]^{R\Ph_h}\ar[u]_{\ttt \OO(-2)}   
                      &\D^b(Q_h)\ar[u]_{RC_h^+}
}
$$
where $RC_h^+$ is the Coxeter element \eqref{e:C+}, and similarly for
$LC_h^-$. 
\end{corollary}

\begin{remark}
It should be noted that while $R\Ph_h$ is an equivalence of derived
categories, it is definitely not true that $\Ph_h$ is an equivalence
of abelian categories. For example, there are no injective or
projective objects in $\C$ while there are enough injectives and
projectives in $\Rep Q_h$. Similarly, the set of simple modules is
very much different in these two categories. However, the sets of
indecomposable modules are effectively the same. 
\end{remark}

\section{The root system and the Grothendieck group}
In this section, we list some important corollaries of the equivalence of
categories constructed in the previous section. Again, many of these
results can also be obtained form the equivalence of categories
constructed by Lenzing; however, we choose to provide an independent
exposition.

Thoroughout this section, we let $L=K(\C)$ be the Grothendieck group of
category $\C$. We define the inner product on $L$ by 
\begin{equation}\label{e:euler_form}
    (x,y)=\<x,y\>+\<y,x\>
\end{equation}
 where by definition 
$$
    \<[X],[Y]\>=\dim \RHom (X,Y)=\dim \Hom (X,Y) -\dim \Ext^1(X,Y)
$$
(compare with \thref{t:reps_of_quivers}). We define $\Delta\subset L$ by 
\begin{equation}
\Delta=\{[X], X \text{--- a non-zero indecomposable object in
}\D^b(\C)\}
\end{equation}

Finally, we define the map $C\colon \Delta\to\Delta$ by 
$$
C([X])=[X(-2)]
$$

\begin{theorem}\label{t:main3}
\par\indent
\begin{enumerate}
  \item The set $\Delta\subset L$ is an affine root system, and $C$ is
    a Coxeter element. 
  \item Recall the lattice $L(Q)$ and root system $\Delta(Q)$ from
    \seref{s:basic}. Then for any height function $h$ on $Q$, the map
    \begin{equation}
    \begin{aligned}
    R\Ph_h\colon L&\to L(Q)\\
              [\F]&\mapsto \bigoplus 
              \biggl(\dim \RHom_\C(X_i(h_i), \F)\biggr) \al_i
    \end{aligned}
    \end{equation}
    is an isomorphism of abelian groups which identifies
    $\Delta\subset L$  with $\Delta(Q)\subset L(Q)$ and $C$ with the
    Coxeter element  $c_h^+$ defined by \eqref{e:coxeter}. 
\end{enumerate}
\end{theorem}
\begin{proof}
  The first part follows from the second one; the
  second part follows from the fact that $R\Ph_h$ is an equivalence of
  categories (\thref{t:main2} and \coref{c:coxeter}). 
\end{proof}

\begin{corollary}\label{c:generators_of_K(C)}
  \par\noindent
  \begin{enumerate}
  \item 
    For any height function $h$, the classes $[X_q]$, $q\in Q_h$, are free
      generators of   the Grothendieck group $K(\C)$. 
   \item 
    $K(\C)$ is generated by the classes $[X_i(n)], n\equiv p(i)\mod 2$,
    and AR relations \eqref{e:aus-reiten3} is the full list of all
    relations among classes $[X_i(n)]$.
  \end{enumerate}
\end{corollary}
\begin{proof}
By \leref{l:generators}, classes $[X_q]$, $q\in Q_h$, generate $K(\C)$. On
the other hand, by \thref{t:main3}, $K(\C)$ is isomorphic to $\Z^I$ and
thus has rank $|I|$, so these generators must be linearly independent.  

To prove the second part, denote temporarily $L'=Span([X_i(n)])/J$, where
$J$ is the subgroup generated by AR relations \eqref{e:aus-reiten3}.
Since $[X_i(n)]$ generate $K(\C)$ and AR relations hold in $K(\C)$, we
see that $K(\C)$ is a quotient of $L'$. 

Now choose some height function $h$. It follows from 
\leref{l:reflection_functors1} that $[X_q]$, $q\in Q_h$, generate $L'$, so
it has rank at most $|I|$. On the other hand, $K(\C)$ has rank $|I|$.
Thus, $L'$ has rank $|I|$ and $L'=K(\C)$. 
\end{proof}

\begin{example}\label{x:standard2}
  Let $h$ be the ``standard'' height function as defined in
\exref{x:standard}.  Then the map
$R\Ph_h\colon K(\C)\to  L(Q)\simeq K(G)$
  is given by 
  \begin{align*}
  &[X_i]\mapsto \al_i, \quad i \in I_0\\
  &[X_i(1)]\mapsto  \al_i+\sum_{j-i}\al_j, \quad i \in I_1
  \end{align*}
  and the corresponding Coxeter element is $C=(\prod_{i\in
    I_1}s_i)(\prod_{i\in I_0}s_i)$. 
  
  This also implies that for this $h$, $[X_i(-1)]\mapsto -\al_i$, $i\in
    I_1$; in particular, we see that classes
  $\al_i, i\in I_0$, and $-\al_i, i\in I_1$, form a set of representatives
  of $C$-orbits on on $\Qhat$. An analogous statement for finite root
  system has been proved in \cite{kostant2}. 
\end{example}

For completeness, we also describe here the indecomposable sheaves
corresponding to imaginary roots of $\Delta$.

\begin{theorem}\label{t:imaginary_root}
Let $x\in \PP$ be a generic point: $\Stab_\Gbar x=\{1\}$, and let 
$$
\de=[\OO_{[Gx]}]\in \Delta
$$
\textup{(}see \eqref{e:torsion_sheaf}\textup{)}. Then: 
\begin{enumerate}
\item $\de$ does not depend on the choice of point $x$
\item $\de=\de_0-\de_1$, where 
\begin{align*}
\de_0&=\sum_{i\in I_0}d_i[X_i]=[R_0]\\
\de_1&=\sum_{i\in I_1}d_i[X_i(-1)]=[R_1(-1)],
\end{align*}
 where $d_i=\dim X_i$ and $R_0, R_1$ are even and odd parts of the
regular representation defined by \eqref{e:R_p}. 

\item $C\de=\de$; $C\de_0=\de_0-2\de$; $C\de_1=\de_1-2\de$
\item For any $\al\in L$, $(\de,\al)=0$
\item $\de$ is a generator of the set of imaginary roots of $\Delta$:
$$
\Delta^{im}=\Z\de
$$
\end{enumerate}
\end{theorem}
\begin{proof}
We start by proving (2), by explicitly constructing a resolution of 
$\OO_{[Gx]}$. 

Let $\ph_x$ be a holomorphic section of $\OO(1)$ which has a simple
zero at $x$. Then we have a short exact sequence of sheaves (not
equivariant):
$$
0\to\OO\xxto{\ph_x}\OO(1)\to \OO_{[x]}\to 0
$$
Tensoring it with $\OO(-1)$, we get a sequence
$$
0\to\OO(-1)\xxto{\ph_x}\OO\to \OO_{[x]}\to 0
$$
Now let us apply to this sequence functor $\Ind$ defined by
\eqref{e:ind}. Using \leref{l:frobenius}, we see that it gives a
$\Gbar$-equivariant short exact sequence 
$$
0\to R_1\ttt \OO(-1)\to R_0\ttt\OO\to \OO_{[Gx]}\to 0
$$
which gives equality $\de=\de_0-\de_1$, thus proving part (2) of the
theorem.

Part (1) follows from (2).

Since $\OO_{[Gx]}\ttt \OO(-2)\simeq \OO_{[Gx]}$, we get $C\de=\de$.
To compute $C\de_0$, note that the same argument as in the proof of
part (2) also gives a short exact sequence 
$$
0\to R_0\ttt \OO\to R_1\ttt\OO(1)\to \OO_{[Gx]}\to 0
$$
thus giving $C^{-1}\de_1=\de+\de_0=2\de+\de_1$. 

To prove part (4), recall notation
$\<[X],[Y]\>=\dim\Hom(X,Y)-\dim\Ext^1(X,Y)$. Then Serre duality
immediately gives 
$$
\<x,y\>=-\<y,Cx\>.
$$
Since $C\de=\de$, we get 
$$
(\de,x)=\<x,\de\>+\<\de,x\>=\<x,\de\>-\<x,C\de\>=0
$$

Part (4) implies that $\de$ is an imaginary root. Moreover, since by
\coref{c:generators_of_K(C)} classes $[X_i]$, $[X_i(-1)]$ are free
generators of $L$. Since some of $d_i$ are equal to 1,  it follows
from part (2) that for any $k>1$, $\de/k\notin L$; thus, $\de$ must be
a generator of the set of imaginary roots. 
\end{proof}

Since every indecomposable object in $\D^b(\C)$ is of the form
$\F[n]$, where $\F$ is an indecomposable $G$-equivariant sheaf
(\thref{t:indecomposable}), it follows that every root $\al\in\De$
can be written as either $\al=[\F]$ or $\al=-[\F]$, thus giving some
splitting of $\De$ into positive and negative roots. This polarization
can be described explicitly. 

Recall from algebraic geometry (see, e.g., \cite[Exercises II.6.10,
II.6.12]{hartshorne}) that for any coherent sheaf, we can
define two integer numbers, its {\em rank} and {\em degree}. In
particular, for a locally free sheaf $\F=X\ttt \OO(n)$ (where $X$ is a
finite-dimensional vector space), we have 
\begin{align*}
&\rk (X\ttt \OO(n))=\dim X,\\
&\deg (X\ttt \OO(n))= n\dim X.
\end{align*}
Degree and rank give  well-defined linear maps $K\to \Z$, where $K$
is the  Grothendieck group of the category of coherent sheaves. 

In particular, we can define rank and degree for a $G$-equivariant
sheaf, ignoring the equivariant structure; this gives linear maps
$K(\C)\to\Z$, which we also denote by $\rk, \deg$.

\begin{lemma}\label{l:rank}
\par\indent
  \begin{enumerate}
    \item If $\F\in\C$ is a non-zero free sheaf, then $\rk \F>0$.
    \item If $\F\in\C$ is a non-zero torsion sheaf, then $\rk \F=0$,
          $\deg\F>0$.  
    \item For any $x\in K(\C)$, we have 
    $$
    \rk (x)=(x,\de_0)=(x,\de_1)
    $$
    where $\de_0$, $\de_1$ are defined in \thref{t:imaginary_root}.
  \end{enumerate}
\end{lemma}
\begin{proof}
  The first two parts are well-known. 
   
  To check  $\rk (x)=(x,\de_0)=(x,\de_1)$, it suffices to check it
  for  $x=[X_i]$, $i\in I_0$ and $x=[X_i(-1)], i\in I_1$ (by
  \exref{x:standard2},  they generate $L$). If $i\in I_0$,
  then  
  $$
  ([X_i], \de_0)=\dim \Hom_G(X_i,\sum_{j\in I_0} d_jX_j)=d_i=\rk X_i
  $$
  thus proving $\rk (x)=(x,\de_0)$. Since $\de_0-\de_1=\de$ is in the
  kernel of $(\, , \,)$ (\thref{t:imaginary_root}), this implies
  $(x,\de_0)=(x,\de_1)$. 
    
  If $i\in I_1$, then 
 $$
  ([X_i(-1)], \de_1)=\dim \Hom_G(X_i(-1),\sum_{j\in I_1}
  d_jX_j(-1))=d_i=\rk X_i
 $$
  thus proving $\rk (x)=(x,\de_1)$. Since $\de_0-\de_1=\de$ is in the
  kernel of $(\, , \,)$ (\thref{t:imaginary_root}), this implies
  $(x,\de_0)=(x,\de_1)$.

\end{proof}

Note that the functional $\deg x$ can not be written in terms of the
form $(\, , \,)$: indeed, $\deg \de =|\Gbar|=|G|/2$, but $(\de,
\cdot)=0$. 

\begin{theorem}
Let $\al\in \De$.
\begin{enumerate}
  \item $\al=[\F]$ for some indecomposable free sheaf $\F\in \C$ 
    iff $\rk(\al)=(\al,\de_0)>0$
  \item $\al=-[\F]$ for some indecomposable free sheaf $\F\in \C$ 
    iff $\rk(\al)=(\al,\de_0)<0$
  \item $\al=[\F]$ for some indecomposable torsion sheaf $\F\in \C$ 
    iff $\rk(\al)=(\al,\de_0)=0$, $\deg(\al)>0$
  \item $\al=-[\F]$ for some indecomposable torsion sheaf $\F\in \C$ 
    iff $\rk(\al)=(\al,\de_0)=0$, $\deg(\al)<0$
\end{enumerate}
\end{theorem}
Thus, we see that we have a triangular decomposition of $\De$:
\begin{equation}\label{e:polarization}
  \begin{aligned}
  &\De=\De_+\sqcup\De_0\sqcup \De_-\\
  &\De_+=\{\al\in \De\st \rk(\al)>0\}=\{[\F],\ \F\text{ --- a free
    indecomposable sheaf}\}\\
  &\quad \text{(note that $\De_+$ is exactly the set of vertices of
    $\Qhat$)}\\
  &\De_-=\{\al\in \De\st \rk(\al)<0\}=\{-[\F],\ \F\text{ --- a free
    indecomposable sheaf}\}\\
  &\De_0=\{\al\in \De\st \rk(\al)=0\}=\{\pm[\F],\ \F\text{ --- a free
    torsion sheaf}\}
  \end{aligned}
\end{equation}
The set $\De_+$ has been discussed by Schiffmann \cite{schiffmann2}, who
used notation $\Qhat_+$ and denoted the corresponding subalgebra in the
loop algebra by $\mathcal{L}\mathfrak{n}$. Note, however, that this
notation is somewhat misleading:  $\mathcal{L}\mathfrak{n}$ is not the
loop algebra of a positive part of the finite dimensional algebra
$\bar\g$, which easily follows from the fact that there are real roots in 
$\De_0$ (see example in the next section).
 
\begin{theorem}\label{t:C^g}
Let $g=|\Gbar|=|G|/2$. Then for any $x\in L$ we have 
\begin{equation}\label{e:C^g}
C^g (x)=x-2(\rk x)\de 
\end{equation}
\end{theorem}

\begin{proof}
Let $\ph_x$ be a section of the sheaf $\OO(1)$ which has a single zero at
generic point $x$. Then we have a short exact sequence 
$$
0\to \OO\xxto{\ph_x} \OO(1)\to \OO_{[x]}\to 0
$$
which, however, is not equivariant even under the action of $\{\pm
I\}\subset \SU$. To fix it we consider $\ph_x^2$ which gives the following
$\Z_2$-equivariant sequence
$$
0\to \OO(-2)\xxto{\ph^2_x} \OO\to \OO_{2[x]}\to 0
$$
Now let us take product of pullbacks of $\ph_x^2$ under all $g\in\Gbar$
$$
0\to \OO(-2g)\xxto{\prod_{g\in\Gbar}g^*\ph^2_x} \OO\to \OO_{2[Gx]}\to 0
$$

Tensoring it with any locally free sheaf $\F$, we get 
$$
0\to \F(-2g)\xxto{\prod_{g\in\Gbar}g^*\ph^2_x} \F\to
\OO_{2[Gx]}\ttt \F_x\to 0
$$
which implies
$C^g[\F]-[\F]+2(\rk \F)\de=0$.
\end{proof}

This result --- that $C^g$ is a translation --- was known before and
can be proved without the use of equivariant sheaves, see e.g.
Steinberg \cite{steinberg}. However, the approach via equivariant
sheaves also provides a nice interpretation for the corresponding
functional as the rank of the sheaf.

Comparing \eqref{e:C^g} with the description of the action of Coxeter
element in the language of representations of the quiver, we see that rank
is closely related to the notion of defect $\partial_c(x)$ as defined in
\cite{drab-ringel}, namely 
$$
\rk(x)=-\frac{1}{2}\partial_c x
$$
Therefore, $\De_0$ is exactly the set of indecomposable objects of defect
zero, which shows that torsion sheaves correspond to regular
representations.

\begin{corollary}
  \par\noindent
  \begin{enumerate}
    \item For any $\al\in \De_0$, $C^gx=x$; in particular, $C$-orbit
      of $\al$ is finite.
    \item $C$ acts freely on $\De_+$, and the set of orbits $\De_+/C$
      is naturally in bijection with $Q$. 
  \end{enumerate}
\end{corollary}

\section{Example: $\widehat{A}_n$}
In this section, we consider the example of the cyclic group of even
order: $G=\Z_n$, $n=2k$. 

The irreducible representations of this group are $X_i$, $i\in \Z_n$;
all of them are one-dimensional. The corresponding Dynkin diagram $Q$
is the cycle of length $n$. 

The root system $\Delta(Q)$ can be described as follows. Let $V$ be a
real vector space of dimension $n+1$, with basis $\delta, e_i$, $i\in
\Z_n$. Define inner product in $V$ by $(e_i,e_j)=\de_{ij}$,
$(v,\de)=0$. Then 
\begin{equation}
  \Delta(Q)=\{e_i-e_j+a\delta, i\ne j\in \Z_n, a\in \frac{j-i}{n}+\Z\}
\end{equation}
The simple roots are 
$$
\al_i=e_i-e_{i+1}+\frac{1}{n}\delta,\quad i\in \Z_n
$$
so that $\sum_i\al_i=\de$. The simple reflections $s_i$ are given by 
\begin{align*}
&s_i(e_i)=e_{i+1}-\frac{1}{n}\delta\\
&s_i(e_{i+1})=e_{i}+\frac{1}{n}\delta\\
&s_i(e_j)=e_j,\quad j\ne i,i+1
\end{align*}

It is easy to see that this description of
$\Delta(Q)$, while unusual, is  equivalent to the standard
description of the affine root system $\widehat{A}_{n-1}$.

We choose standard height function $h$:
\begin{equation}
h(i)=\begin{cases}0,& i\text{ even}\\
                  1, &i\text{ odd}
     \end{cases}
\end{equation}
The corresponding Coxeter  element $C$ is 
$$
C=\left(\prod_{i \text{  odd}} s_i\right)
  \left(\prod_{i \text{ even}} s_i\right)
$$
The action of $C$ on $\Delta(Q)$ is given by 
\begin{equation*}
C(e_i)=\begin{cases}e_{i+2}-\frac{2}{n}\delta,& i\text{ even}\\
                    e_{i-2}+\frac{2}{n}\delta,& i\text{ odd}
       \end{cases}
\end{equation*}
Thus, we have 
\begin{equation*}
C^{n/2}(e_i)=\begin{cases}e_{i}-\delta,& i\text{ even}\\
                    e_{i}+\delta,& i\text{ odd}
       \end{cases}
\end{equation*}
which implies
\begin{equation}
\begin{aligned}
&C^{n/2}(\al)=\al-(\al,\eps)\delta,\\
&\quad \eps = \sum_{i \text{ even}}e_i-\sum_{i \text{ odd}}e_i\equiv 
\sum_{i \text{ even}}\al_i\equiv -\sum_{i \text{ odd}}\al_i\mod \Z\delta
\end{aligned}
\end{equation}
(compare with \thref{t:C^g}, \leref{l:rank} ). 

Explicitly, we can write 
\begin{align*}
&C^{n/2}(\al)=\begin{cases}
                 \al, &i\equiv j\mod 2 \\
                 \al-2\de, &(i,j)\equiv (0,1)\mod 2 \\
                 \al+2\de, &(i,j)\equiv (1,0)\mod 2 
       \end{cases},\\
&\quad \al=e_i-e_j+a\delta
\end{align*}
Thus, in this case we have 
\begin{align*}
\De_0=\{e_i-e_j+a\de\}, \quad i\equiv j\mod 2\\
\De_+=\{e_i-e_j+a\de\}, \quad i\text{ even}, j\text{ odd}\\
\De_-=\{e_i-e_j+a\de\}, \quad i\text{ odd}, j\text{ even}\\
\end{align*}

\begin{figure}[ht]
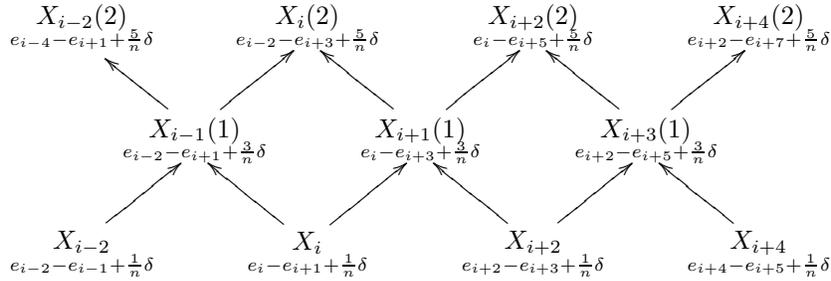

$$ 
\xy
(0,0)*+{X_{i-2}}="A"+(0,-3)*+{\scriptstyle
                      e_{i-2}-e_{i-1}+\frac{1}{n}\delta};
(30,0)*+{X_{i}}="B"+(0,-3)*{\scriptstyle
                            e_{i}-e_{i+1}+\frac{1}{n}\delta};
(60,0)*+{X_{i+2}}="C"+(0,-3)*{\scriptstyle
                           e_{i+2}-e_{i+3}+\frac{1}{n}\delta};
(90,0)*+{X_{i+4}}="D"+(0,-3)*{\scriptstyle
                           e_{i+4}-e_{i+5}+\frac{1}{n}\delta};
(15,15)*+{X_{i-1}(1)}="E"+(0,-3)*{\scriptstyle
                           e_{i-2}-e_{i+1}+\frac{3}{n}\delta}="EE";
(45,15)*+{X_{i+1}(1)}="F"+(0,-3)*{\scriptstyle
                           e_{i}-e_{i+3}+\frac{3}{n}\delta}="FF";
(75,15)*+{X_{i+3}(1)}="G"+(0,-3)*{\scriptstyle
                           e_{i+2}-e_{i+5}+\frac{3}{n}\delta}="GG";
(0,30)*+{X_{i-2}(2)}="H"+(0,-3)*+{\scriptstyle
                      e_{i-4}-e_{i+1}+\frac{5}{n}\delta}="HH";
(30,30)*+{X_{i}(2)}="I"+(0,-3)*{\scriptstyle
                            e_{i-2}-e_{i+3}+\frac{5}{n}\delta}="II";
(60,30)*+{X_{i+2}(2)}="J"+(0,-3)*{\scriptstyle
                           e_{i}-e_{i+5}+\frac{5}{n}\delta}="JJ";
(90,30)*+{X_{i+4}(2)}="K"+(0,-3)*{\scriptstyle
                           e_{i+2}-e_{i+7}+\frac{5}{n}\delta}="KK";
{\ar "A";"EE"};
{\ar "B";"EE"};
{\ar "B";"FF"};
{\ar "C";"FF"};
{\ar "C";"GG"};
{\ar "D";"GG"};
{\ar "E";"HH"};
{\ar "E";"II"};
{\ar "F";"II"};
{\ar "F";"JJ"};
{\ar "G";"JJ"};
{\ar "G";"KK"};
\endxy
$$
 \caption{Fragment of graph $\Qhat$ for root system of type $A$.
 Here $i$ is an even number. The figure also shows, for each $X_q\in
 \Qhat$, the image of  $R\Ph_h[X_q]$ for the standard choice of height
 function $h$ as in \exref{x:standard}.} \label{f:An_graph}
\end{figure}

\firef{f:An_graph} shows a segment of the AR graph $\Qhat$ for
this root system.

\bibliographystyle{amsalpha}


\end{document}